\documentclass[11pt]{article}
\usepackage[a4paper]{anysize}\marginsize{3.5cm}{3.5cm}{1.3cm}{2cm}
\pdfpagewidth=\paperwidth \pdfpageheight=\paperheight
\usepackage{amsfonts,amssymb,amsthm,amsmath,eucal}
\usepackage{pgf}
\usepackage{bbm,array}
\theoremstyle{plain}
\newtheorem{thm}{Theorem}[section]
\newtheorem{theorem}[thm]{Theorem}

\newtheorem{lemma}[thm]{Lemma}

\newtheorem{conjecture}[thm]{Conjecture}

\theoremstyle{definition}

\newtheorem{question}[thm]{Question}

\newtheorem{thevarthm}[thm]{\varthmname}

\newenvironment{varthm*}[1]{\trivlist\item[]{\bf #1.}\it}{\endtrivlist}

\renewcommand\geq{\geqslant}

\renewcommand\leq{\leqslant}

\newcommand\be{\begin{eqnarray*}}
\newcommand\ee{\end{eqnarray*}}

\newcommand\Q{\mathbb Q}

\newcommand\C{\mathbb C}

\renewcommand\P{\mathbb P}

\newcommand\cali{{\mathcal I}}
\newcommand\calj{{\mathcal J}}
\newcommand\calk{{\mathcal K}}

\def\field{\C}
\newcommand\newop[2]{\def#1{\mathop{\rm #2}\nolimits}}
\newop\log{log}
\newop\ord{ord}
\newop\Gal{Gal}
\newop\SL{SL}
\newop\Bl{Bl}
\newop\mult{mult}
\newop\mass{mass}
\newop\div{div}
\newop\codim{codim}
\newop\sing{sing}
\newop\vdim{vdim}
\newop\edim{edim}
\newop\Ass{Ass}
\newop\size{size}
\newop\reg{reg}
\newop\satdeg{satdeg}
\newop\supp{supp}
\newcommand\eqnref[1]{(\ref{#1})}

\newcommand\eps{\varepsilon}
\newcommand\lra{\longrightarrow}

\def\keywordname{{\bfseries Keywords}}%
\def\keywords#1{\par\addvspace\medskipamount{\rightskip=0pt plus1cm
\def\and{\ifhmode\unskip\nobreak\fi\ $\cdot$
}\noindent\keywordname\enspace\ignorespaces#1\par}}
\def\subclassname{{\bfseries Mathematics Subject Classification
(2000)}\enspace}
\def\subclass#1{\par\addvspace\medskipamount{\rightskip=0pt plus1cm
\def\and{\ifhmode\unskip\nobreak\fi\ $\cdot$
}\noindent\subclassname\ignorespaces#1\par}}

\begin{document}

\author{M.~Dumnicki, T.~Szemberg\footnote{The second named author was partially supported
by NCN grant UMO-2011/01/B/ST1/04875}, H.~Tutaj-Gasi\'nska}
\title{Counterexamples to the $I^{(3)}\subset I^2$ containment}
\date{\today}
\maketitle
\thispagestyle{empty}

\begin{abstract}
   The purpose of this short note is to
   show that there is in general no containment
   $$I^{(3)}\subset I^2$$
   for an ideal $I$ of points in $\P^2$.
   This answers in the negative a question asked by Huneke
   and generalized by Harbourne.
   The sets of points constituting  counter-examples
   come from the dual of the Hesse configuration
   and more generally from Fermat arrangements.
\subclass{MSC 14C20 \and MSC 13C05 \and MSC 14N05 \and MSC 14H20 \and MSC 14A05}
\end{abstract}


\section{Introduction}
\label{intro}
   Let
   $\cali\subset S=\C[x_0,\dots,x_n]$ be a homogeneous ideal
   in the graded ring of polynomials.
   The $m$--th \emph{symbolic power} $\cali^{(m)}$ of $\cali$
   is defined as
   $$
   \cali^{(m)} = S \cap \left( \bigcap_{\mathfrak{p} \in \Ass(\cali)} \cali^{m} S_{\mathfrak{p}} \right),
   $$
   where the intersection is taken in the field of fractions of $S$.

   There has been considerable interest in containment relations
   between usual and symbolic powers of homogeneous ideals over the
   last two decades. The most general results in this direction
   have been obtained with multiplier ideal techniques in characteristic
   zero by Ein, Lazarsfeld and Smith \cite{ELS01} and using tight
   closures in positive characteristic by Hochster and Huneke \cite{HoHu02}.
   Applying these results to a homogeneous ideal $\cali$ in the coordinate
   ring $S$ of the projective space we obtain the following
   containment statement
\begin{equation*}\label{eq:els containment}
   \cali^{(nr)}\subset \cali^r\;\mbox{ for all } r\geq 0.
\end{equation*}
   Quite a number of examples has suggested that the
   following statement could be true, see \cite[Conjecture 8.4.2]{PSC},
   \cite[Conjecture 1.1]{BocHar10b}, \cite[Conjecture 4.1.1]{HaHu}.
\begin{conjecture}\label{conj:hahu}
   Let $\cali\subset\C[\P^n]$ be a homogeneous ideal.
   For $m\geq rn-(n-1)$ there is the containment
   $$\cali^{(m)}\subset \cali^r.$$
\end{conjecture}
   This conjecture asserts in particular that an earlier question raised by Huneke
   has a positive answer; see \cite[Problem 0.4]{Hun06}, see also \cite[page 400]{BocHar10a}
   and \cite[section 4.1]{HaHu}. The referee has kindly informed us that Huneke had
   been raising this question verbally for a few years before 2006 but \cite[Problem 0.4]{Hun06}
   seems to be its first occurrence in print.
\begin{question}[Huneke]\label{huneke}
   Let $\cali$ be a homogeneous radical ideal of points in the projective plane.
   Is there then the containment
   $$\cali^{(3)}\subset \cali^2\,?$$
\end{question}
   This question has been affirmatively answered for general points, see \cite{BocHar10a},
   star configurations, see \cite{HaHu}, complete intersections and some special configurations
   of points, see \cite{BCH}.

   We show here that the containment in Question \ref{huneke} fails in general.
   This implies also that Conjecture \ref{conj:hahu} is false. There is quite
   a number of closely related yet distinct conjectures concerning containment relations between symbolic
   and usual powers of ideals in the literature. For the convenience of the reader
   we mention here that the results of this note show that
   the following conjectures are false:
   \cite[Conjecture 4.1.1 and Conjecture 4.1.5]{HaHu} (for $N=r=2$),
   \cite[Conjecture 8.4.2]{PSC}, \cite[Conjecture 1.1]{BocHar10b},
   \cite[Conjecture 3.9]{BCH} (for $N=t=2$ and $m=1$),
   also Questions 4.2.2 and 4.2.3 in \cite{BCH} have a negative answer.

\section{The dual Hesse configuration}
\label{sec:Hesse}
   We begin with an explicit realization of the dual Hesse
   configuration. Note that up to projective change of coordinates
   there is a unique configuration of that type \cite[Example 7.3]{Pal11}.

   Let $\eps$ be a primitive root of $1$ of order $3$.
   We consider the radical ideal $\cali$
   of the following set of $12$ points in $\P^2$:
   \begin{equation*}
   \begin{array}{lll}
      P_1=(1:0:0),         & P_2=(0:1:0),         & P_3=(0:0:1),\\
      P_4=(1:1:1),         & P_5=(1:\eps:\eps^2), & P_6=(1:\eps^2:\eps),\\
      P_7=(\eps:1:1),      & P_8=(1:\eps:1),      & P_9=(1:1:\eps),\\
      P_{10}=(\eps^2:1:1), & P_{11}=(1:\eps^2:1), & P_{12}=(1:1:\eps^2).
   \end{array}
   \end{equation*}
   These points form a $12_39_4$ configuration, i.e. there are $9$
   lines
   \begin{equation*}
   \begin{array}{lll}
   L_1:\; x-y,        & L_2:\; y-z,        & L_3:\; z-x,\\
   L_4:\; x-\eps y,   & L_5:\; y-\eps z,   & L_6:\; z-\eps x,   \\
   L_7:\; x-\eps^2 y, & L_8:\; y-\eps^2 z, & L_9:\; z-\eps^2 x.
   \end{array}
   \end{equation*}
   such that exactly $3$ configuration lines pass through each of configuration points
   and exactly $4$ points lie on a configuration line. This is the dual
   of the well known Hesse configuration, see \cite{ArtDol09}
   for a lot more on this beautiful subject.

   Turning back to the ideal $\cali$, we exhibit first its
   generators.
\begin{lemma}\label{lem:generators}
   The ideal $\cali$ is generated by polynomials
   $$f_1 := z(x^3-y^3),\; f_2 := x(y^3-z^3)\;\mbox{ and }\; f_3=y(z^3-x^3).$$
\end{lemma}
\proof
   We have obviously $(f_1,f_2,f_3)\subset\cali$, so it remains to check
   the opposite inclusion. To this end let first
   $\calj$ be the radical ideal of points
   $P_4,\dots,P_{12}$.\\
   \textbf{Claim.} Polynomials
   $$g_1=z^3-x^3,\;\mbox{ and }\; g_2=y^3-z^3$$
   generate $\calj$.\\
   Let $g \in \calj$ be a homogeneous element. Using identities
   $$x^3 = -(z^3-x^3) + z^3, \quad y^3 = (y^3-z^3) + z^3$$
   we can write $g=g'+g''$ with some homogeneous $g'$, $g''$ satisfying
   $g'' \in (g_1,g_2)$ and $g'$ such that $\deg_{x} g' \leq 2$ and
   $\deg_{y} g' \leq 2$. Note that $g'=g-g''\in\calj$.

   Now we set $h(x,y)=g'(x,y,1)$ and let $\calk$ be the radical ideal of points $P_1,\dots,P_{12}$
   in the affine chart $z=1$. Then $h\in\calk$. Note that $h$ is supported on the following set of monomials
   $$\supp(h) \subset Q=\{1,x,y,x^2,xy,y^2,x^2y,xy^2,x^2y^2\}.$$
   Observe that the set $Q$ is a monomial basis for the algebra $\field[x,y]/\calk$.
   Indeed, since $(x-1)(x-\eps)(x-\eps^2)$ and $(y-1)(y-\eps)(y-\eps^2)$ belong to $\calk$
   the set $Q$ generates $\field[x,y]/\calk$. Since
   $\dim_{\field} \field[x,y]/\calk = 9$ (we have 9 points), $Q$ is a basis.

   Now the inclusion $h \in \calk$ with $h$ supported on $Q$ implies $h=0$, thus
   $g'=0$ and consequently $g=g'' \in (g_1,g_2)$. This exactly means that $\calj = (g_1,g_2)$.

   Now we turn back to the inclusion
   $$\cali\subset (f_1,f_2,f_3).$$
   All polynomials in the subsequent argument are supposed to be homogeneous.
   Let $g \in \cali$ be an arbitrary element.
   Of course $g \in \calj$, hence
   $$g = h_1(x,y,z)\cdot (z^3-x^3) + h_2(x,y,z)\cdot (y^3-z^3)$$
   for some polynomials $h_1$ and $h_2$.
   We can split
   $$h_1(x,y,z) = zh_3(x,y,z) + h_4(x,y)$$
   into monomials containing $z$ and those depending only on $x$ and $y$.
   Hence we can also write
   $$g = h_4(x,y)\cdot (z^3-x^3) + (h_2(x,y,z)-zh_3(x,y,z))\cdot (y^3-z^3) - h_3(x,y,z)\cdot z(x^3-y^3).$$
   Gathering together terms divisible by $y(z^3-x^3)$, $x(y^3-z^3)$ and $z(x^3-y^3)$ we can write
   $$g = h_5(x)(z^3-x^3) + h_6(y,z)(y^3-z^3) + h_7(x,y,z)$$
   for some $h_7 \in (f_1,f_2,f_3)$. Obviously
   $$h_5(x)(z^3-x^3) + h_6(y,z)(y^3-z^3)$$
   vanishes at $P_1$, $P_2$, $P_3$. This implies
   $$h_5(x) = 0\;\mbox{ and consequently }\; yz\;\mbox{ divides }\;h_6(y,z).$$
   Since
   $$yz(y^3-z^3) = -z \cdot y(z^3-x^3) - y \cdot z(x^3-y^3) \in (f_1,f_2,f_3)$$
   this completes the proof.
\endproof
   We have the following two relations between the generators of $\cali$:
   $$xyf_1+yzf_2+zxf_3=0\;\mbox{ and }\; z^2f_1+x^2f_2+y^2f_3=0.$$
   It is easy to check that these relations determine the minimal resolution of $\cali$:
   $$0\lra\bigoplus\limits^2 S(-6)\lra\bigoplus\limits^3 S(-4)\lra \cali\lra 0.$$
   Hence 
   the Castelnuovo-Mumford regularity of $\cali$
   is $\reg(\cali)=5$. Then
   \cite[Theorem 1.1]{GGP95} implies that
   $$\left(\cali^{(2)}\right)_t = (\cali^2)_t\;\mbox{ for }\; t\geq 10$$
   and hence
   $$\left(\cali^{(3)}\right)_t\subset (\cali^{(2)})_t= (\cali^2)_t \mbox{ for }\; t\geq 10.$$
   Thus the containment problem $\cali^{(3)}\subset\cali^2$
   reduces to finding an element of degree less than $10$ in $\cali^{(3)}$ which is not in $I^2$.
\begin{theorem}\label{thm:containment}
   The polynomial
   $$f=L_1\cdot\ldots\cdot L_9=x^3y^6-x^6y^3+y^3z^6-y^6z^3+z^3x^6-z^6x^3=(x^3-y^3)(y^3-z^3)(z^3-x^3)$$
   is contained in $\cali^{(3)}$ but it is not contained in
   $\cali^2$.
\end{theorem}
\begin{proof}
   The geometry of the configuration implies that
   $f\in \cali^{(3)}$. In fact this is the only nonzero element (up to a multiplicative constant)
   of degree $9$ in the third symbolic power of $\cali$, which can be easily checked by Bezout's theorem.

   For the second claim we assume to the contrary that $f \in \cali^2$.
   Using Lemma \ref{lem:generators} we can write
   \begin{equation}\label{fini}
   f = (ax+by+cz)f_1^2 + \text{ other terms divisible by } xy \text{ or } z^4.
   \end{equation}
   Substituting $x=0$ in \eqnref{fini} we get
   $$y^3z^6-y^6z^3 = (by+cz)y^6z^2 + \text{terms divisible by } z^4.$$
   Comparing coefficients at $y^6z^3$ we obtain $c=-1$.

   Substituting in turn $y=0$ in \eqnref{fini} we get
   $$z^3x^6-z^6x^3 = (ax+cz)x^6z^2 + \text{terms divisible by } z^4,$$
   which comparing again coefficients at $x^6z^3$ gives $c=1$, a contradiction.
\end{proof}
   One can use the following Singular \cite{DGPS} script
   in order to verify all above claims, in particular that
   $f\notin\cali^2$.
   \begin{verbatim}
   ring R=(0,e),(x,y,z),dp; option(redSB);
   minpoly=e2+e+1;
   ideal P1=y,z; ideal P2=x,z; ideal P3=x,y;
   ideal P4=x-z,y-z; ideal P5=y-e*x,z-e^2*x; ideal P6=y-e^2*x,z-e*x;
   ideal P7=x-e*z,y-z; ideal P8=x-z,y-e*z; ideal P9=y-x,z-e*x;
   ideal P10=x-e^2*z,y-z; ideal P11=x-z,y-e^2*z; ideal P12=y-x,z-e^2*x;
   ideal I=intersect(P1,P2,P3,P4,P5,P6,P7,P8,P9,P10,P11,P12);
   regularity(mres(I,0)); 
   ideal I3=intersect(P1^3,P2^3,P3^3,P4^3,P5^3,
                      P6^3,P7^3,P8^3,P9^3,P10^3,P11^3,P12^3);
   poly F=I3[1];
   reduce(F,std(I^2));
   quit;
   \end{verbatim}

\section{Fermat arrangements}
\label{sec:Fermat}
   The dual Hesse configuration from section \ref{sec:Hesse} is a special
   case of Fermat arrangements, see \cite[Example II.6]{Urz08}.
   More specifically, the $3d$ lines in a Fermat arrangement are defined
   as the zero locus of the polynomial
   $$f_d=(x^d-y^d)(y^d-z^d)(z^d-x^d).$$
   If $\eta$ is a primitive root of $1$ of order $d$,
   then the intersection points of these lines are
   $Q_{a,b}=(1:\eta^a:\eta^b)$ for $a,b=1\dots,d$
   and the three coordinate points $P_1=(1:0:0)$, $P_2=(0:1:0)$
   and $P_3=(0:0:1)$. There are exactly $d$ lines meeting
   in a coordinate point and exactly $3$ lines passing through
   every point $Q_{a,b}$. Taking $\cali_d$ as the radical ideal
   of the union of all points $Q_{a,b}$ and $P_i$, it follows
   that $f_d\in\cali^{(3)}$. It turns out that
   $f_d$ is not an element of $\cali^2$. So also in this case
   we have
   $$\cali^{(3)}\not\subset\cali^2.$$
   The arguments are
   similar as in section \ref{sec:Hesse} and we don't pursue
   any exact proofs here. Note that for $d=3$ we recover exactly the dual
   Hesse configuration.

   It is natural to wonder if Conjecture \ref{conj:hahu}
   could be true after some modifications. The constructions carried out
   in this note do not exclude the following variant of Conjecture \ref{conj:hahu}.
\begin{question}
   Let $\cali$ be a homogeneous radical ideal of points in the projective
   plane. Is there then the containment
   $$\cali^{(m)}\subset \cali^r$$
   for $m\geq 2r-1$ and $r\geq 3$?
\end{question}
   The first problem to decide would be the containment $\cali^{(5)}\subset\cali^3$.

\paragraph*{\emph{Acknowledgement.}}
   We would like to thank Brian Harbourne for helpful remarks.
   We thank also the referee for numerous remarks which improved
   the readability of this note.



\bigskip \small

\bigskip
   Marcin Dumnicki,
   Jagiellonian University, Institute of Mathematics, {\L}ojasiewicza 6, PL-30-348 Krak\'ow, Poland

\nopagebreak
   \textit{E-mail address:} \texttt{Marcin.Dumnicki@im.uj.edu.pl}

\bigskip
   Tomasz Szemberg,
   Instytut Matematyki UP,
   Podchor\c a\.zych 2,
   PL-30-084 Krak\'ow, Poland.

\nopagebreak
   \textit{E-mail address:} \texttt{szemberg@up.krakow.pl}

\bigskip
   Halszka Tutaj-Gasi\'nska,
   Jagiellonian University, Institute of Mathematics, {\L}ojasiewicza 6, PL-30-348 Krak\'ow, Poland

\nopagebreak
   \textit{E-mail address:} \texttt{Halszka.Tutaj@im.uj.edu.pl}


\end{document}